\documentclass[12pt,reqno]{amsart}
\usepackage{amsmath,amssymb,amsfonts,amsthm}
\usepackage[mathscr]{eucal}
\usepackage{diagrams}
\usepackage{hyperref}

\author{Theodore Voronov}

\address{Department of Mathematics, University of Manchester Institute of Science and Technology
(UMIST), United Kingdom.} \address{Phone: +44 161 200 3682, fax: +44
161 200 3669}

\email{theodore.voronov@umist.ac.uk}

\title[Higher derived brackets]{Higher derived brackets and homotopy algebras}

\newtheorem{thm}{Theorem}

\newtheorem{cor}{Corollary}
\newtheorem{lem}{Lemma}

\theoremstyle{definition}
\newtheorem{de}{Definition}
\theoremstyle{remark}
\newtheorem{ex}{Example}[section]

\newtheorem{rem}{Remark}[section]
\newtheorem*{Rem}{Remark}

\def\co{\colon\thinspace}

\renewcommand{\leq}{\leqslant}
\renewcommand{\geq}{\geqslant}

 \DeclareMathOperator{\ord}{ord}
 
 \DeclareMathOperator{\End}{End}
 \DeclareMathOperator{\Vect}{Vect}
\DeclareMathOperator{\ad}{ad} \DeclareMathOperator{\Ker}{Ker}
\renewcommand{\Im}{{\mathop{\mathrm{Im}}}}
\newcommand{\void}{\varnothing}
\newcommand{\V}{{\mathfrak{V}}}
\newcommand{\id}{{\mathrm{id}}}

\newcommand{\lie}[1]{{\mathcal L}_{{#1}}}
\newcommand{\der}[2]{{\frac{\partial {#1}}{\partial {#2}}}}

\newcommand{\Z}{{\mathbb Z_{2}}}
\newcommand{\ZZ}{{\mathbb Z}}
\newcommand{\p}{\partial}
\newcommand{\fun}{C^{\infty}}
\newcommand{\lsch}{{[\![}}
\newcommand{\rsch}{{]\!]}}
\newcommand{\LHS}{{LHS }}

\def\a{\alpha}
\def\e{\varepsilon}
\def\s{\sigma}

\def\D{\Delta}

\newcommand{\F}{{\Phi}}
\renewcommand{\P}{{\Pi}}

\renewcommand{\l}{{\lambda}}
\newcommand{\x}{{\xi}}

\newcommand{\ft}{{\tilde f}}

\newcommand{\itt}{{\tilde \imath}}
\newcommand{\jtt}{{\tilde \jmath}}
\newcommand{\lt}{{\tilde l}}
\newcommand{\gt}{{\tilde g}}
\newcommand{\hht}{{\tilde h}}

\renewcommand{\itt}{{\tilde\imath}}

\newcommand{\at}{{\tilde a}}
\newcommand{\bt}{{\tilde b}}
\newcommand{\ct}{{\tilde c}}

\newcommand{\xt}{{\tilde x}}

\begin{document}

\begin{abstract}
We give a construction of homotopy algebras based on ``higher derived
brackets''. More precisely, the data include a Lie superalgebra with
a projector on an Abelian subalgebra satisfying a certain axiom, and
an odd element $\Delta$.  Given this, we introduce an infinite
sequence of higher brackets on the image of the projector, and
explicitly calculate their Jacobiators in terms of $\Delta^2$. This
allows to control higher Jacobi identities in terms of the ``order''
of $\Delta^2$. Examples include Stasheff's strongly homotopy Lie
algebras and variants of homotopy Batalin--Vilkovisky algebras. There
is a generalization with $\D$ replaced by an arbitrary odd
derivation. We discuss applications and links with other
constructions.
\end{abstract}

\maketitle

\noindent \small{\textbf{Keywords:} Strongly homotopy Lie algebra,
homotopy algebra,  derived bracket,  Batalin--Vilkovisky algebra,
Poisson bracket, homotopy fiber}

\section{Introduction}

Strong  homotopy Lie algebras (``strongly homotopy'', sh Lie
algebras, $L_{\infty}$-algebras) were defined by Lada and Stasheff
in~\cite{lada:stasheff}  (see also~\cite{lada:shla}). According to
Stasheff (private communication), this notion was ``recognized'' by
him when algebraic structures such as  string products  of Zwiebach
(see~\cite{zwiebach:93csft}), and similar, started to appear in
physical works. Before that, Schlessinger and
Stasheff~\cite{schless:stasheff85} realized that the notion of
$L_{\infty}$-algebra was relevant to describing the higher order
obstructions occurring in deformation theory, though this was not
described in the paper~\cite{schless:stasheff85}. Notice also the
work by Retakh~\cite{retakh:liemassey}. The associative counterpart
of the $L_{\infty}$-algebras, Stasheff's $A_{\infty}$-algebras became
widely known much earlier. Currently, all kinds of homotopy algebras
and structures related to them attract great attention. In part, this
is due to their applications such as in Kontsevich's proof of the
existence of deformation quantization for any Poisson manifold. For
an operadic approach to such algebras, see~\cite{mshst:operads}.

In this paper, we give a rather general algebraic construction that
produces strong  homotopy Lie algebras (and related algebras) from
simple data. Namely, we consider a Lie superalgebra $L$ with a
projector on an Abelian subalgebra obeying a ``distributivity''
condition~\eqref{eqdistrib}. There are many examples of such
projectors. Now, given this, an  element $\D$ defines a sequence of
$n$-ary brackets on the image of the projector $P$ as
\begin{equation*}
\{a_1,\ldots,a_n\}:=  P[\ldots[[\D,a_1],a_2],\ldots,a_n]
\end{equation*}
where $a_i$ are in the image of $P$. We call them \textit{higher
derived brackets} and we call $\Delta$ the \textit{generator} for the
derived brackets.  We prove that for an odd $\D$, the $n$-th
Jacobiator of these derived brackets (i.e., the LHS of the $n$-th
Jacobi identity of the $L_{\infty}$-algebras) exactly equals the
$n$-th derived bracket for the element $\D^2$. Hence, if $\D^2=0$,
our construction leads to strong homotopy Lie algebras. We can weaken
the condition $\D^2=0$ still obtaining the Jacobi identities of
higher orders. This naturally occurs in examples. Particularly
interesting   applications of this construction are to higher Poisson
brackets and  brackets generated by a differential operator, which
give an important  example of a (strong) ``homotopy
Batalin--Vilkovisky algebra''. Our construction  as a particular case
contains the well-known description of $L_{\infty}$-algebras in terms
of homological vector fields. Though it is a generating element $\D$
that  plays a key role in the main examples, it is also possible to
give a similar construction of higher derived brackets taking as a
starting point an arbitrary odd derivation $d\co L\to L$; in
particular, this allows to give a homotopy-theoretic interpretation
of  higher derived brackets.

In Section~\ref{secsetup} we introduce the setup and recall the
notion of $L_{\infty}$-algebras (in a form convenient for our
purposes). In Section~\ref{secmain} we state and prove the main
theorem. Sections~\ref{secsetup} and~\ref{secmain} are purely
algebraic and self-contained.  In Section~\ref{secappl} we consider
some examples of applications. In Section~\ref{secnonin} we return to
algebra, giving a sketch of the generalization of our construction
for non-inner derivations and applying it to  a homotopy-theoretic
interpretation. Finally, in Section~\ref{secdiscuss} we discuss
related works, links with our results and directions for further
study. (Among other things we explain the role of $P$ and the
necessity of higher brackets,  compared to a binary derived bracket
as in~\cite{yvette:derived}.)

\medskip
\textit{Terminology and notation.} We work in the $\Z$-graded (super)
context, e.g., a vector space means a `$\Z$-graded vector space',
etc. Tilde over a symbol denotes parity. (A parallel treatment for
the $\ZZ$-graded context is possible.)

\medskip
\textit{Acknowledgements.} I am deeply grateful to Hovhannes
Khudaverdian, Kirill Mackenzie  and Taras Panov for  stimulating
discussions. Special thanks go to Martin Markl and Jim Stasheff for
their  remarks on the first version of this text. I am particularly
grateful to Jim Stasheff for his detailed comments and for suggesting
numerous improvements of style. An anonymous referee of the paper
suggested to look into a relation between higher derived brackets and
homotopical algebra in the spirit of Quillen. I want to thank him for
this fruitful idea.

\section{Setup and preliminaries} \label{secsetup}

Let $L$ be a Lie superalgebra. Consider a linear projector $P\in
\End L$, $P^2=P$, such that the image of $P$ is an Abelian
subalgebra:
\begin{equation}\label{eqabelian}
    [Pa, Pb]=0
\end{equation}
for all $a,b\in L$. Let $P$ also satisfy the following
distributive law w.r.t. the commutator:
\begin{equation}\label{eqdistrib}
    P[a,b]=P[Pa,b]+P[a,Pb].
\end{equation}
This identity is a convenient way of expressing the requirement that
the kernel of $P$  is also a subalgebra in $L$ (not necessarily
Abelian). Consider an arbitrary odd element $\D$ in $L$. Using these
data, $P$ and $\D$, we shall introduce a sequence of $n$-ary brackets
on the vector space $P(L)\subset L$, the image of the projector $P$,
and check that upon certain conditions they will make it into a
strongly homotopy Lie algebra. More precisely, we shall see how the
corresponding identities are controlled by the properties of the
element $\D^2=\frac{1}{2}[\D,\D]$ and arise step by step.

Let us give some examples of a projector $P$.

\begin{ex}\label{exvectorfields}
Let $V=V_0\oplus V_1$ be a $\Z$-graded vector space, which we also
treat as  a supermanifold. The origin $0$ is a distinguished point.
Take as $L$ the superalgebra $\Vect (V)$ of all vector fields on the
supermanifold $V$ w.r.t. the usual Lie bracket. Let $P$ take every
vector field to its value at the origin considered as a vector field
with constant coefficients. One can check that the map $P\co X\mapsto
X(0)$ satisfies~\eqref{eqdistrib}.
\end{ex}

\begin{ex}\label{exoperators}
Let $A$ be a commutative associative algebra with a unit, and let
$L=\End A$ (the space of all linear operators in $A$) with the usual
commutator of operators as a bracket. The map $P\co \Delta \mapsto
\Delta (1)$ maps every operator to an element of $A$, which can be
identified with an operator of left multiplication. The image of $P$
is an Abelian subalgebra in $\End A$. Again, a direct check shows
that $P$ satisfies~\eqref{eqdistrib}.
\end{ex}

\begin{ex}\label{exhamiltonians}
Let $M$ be a supermanifold, and $T^*M$ its cotangent bundle. Take
as $L$ the Lie superalgebra $\fun (T^*M)$ w.r.t. the canonical
Poisson bracket. Define $P$ as the pullback of functions on $T^*M$
to $M$. $\fun(M)$ can be treated as a subspace of $\fun(T^*M)$; in
particular, it is an Abelian subalgebra. It is directly checked
that $P$ satisfies~\eqref{eqdistrib}. In view of the relation
between the commutator of operators and the Poisson bracket, this
example can be seen as a `classical counterpart' of
Example~\ref{exoperators}.
\end{ex}

Let us recall the definition of a strongly homotopy Lie algebra
due to Stasheff. In a form convenient for our purposes it reads as
follows.

\begin{de}\label{defshL}
A vector space $V=V_0\oplus V_1$ endowed with a sequence of odd
$n$-linear operations, $n=0,1,2,3,\ldots\,$, (which we denote by
braces), is a \textit{\textnormal{(}strongly\textnormal{)}
homotopy Lie algebra} or \textit{$L_{\infty}$-algebra} if: (a) all
operations are symmetric in the $\Z$-graded sense:
\begin{equation}
    \{a_1,\ldots, a_i,a_{i+1}, \ldots, a_n\}=(-1)^{\at_i\at_{i+1}}\{a_1,\ldots, a_{i+1},a_i, \ldots,
    a_n\},
\end{equation}
and (b)  the ``generalized Jacobi identities''
\begin{equation}\label{eqjacn}
    \sum_{\parbox{1.2cm}{\small $\scriptstyle k+l=n$}} \sum_{\text{$(k,l)$-shuffles}}
    (-1)^{\a} \{\{a_{\s(1)},\ldots,a_{\s(k)}\},a_{\s({k+1})},\ldots,a_{\s({k+l})}\}=0
\end{equation}
hold for all $n=0,1,2,\ldots\ $. Here $(-1)^{\a}$ is the sign
prescribed by the sign rule for a permutation of homogeneous
elements $a_1,\ldots, a_n\in V$.
\end{de}

Henceforth symmetric will mean $\Z$-graded symmetric.

The notation is such that the parity of each operation ``sits'' at
the opening bracket, which should be regarded as an odd symbol
w.r.t. the sign rule. A $0$-ary
bracket  is just a distinguished element $\F:= \{\void\}$ in $L$.
Recall that a $(k,l)$-shuffle is a permutation of indices
$1,2,\ldots, k+l$ such that $\s(1)<\ldots <\s(k)$ and
$\s(k+1)<\ldots \s(k+l)$. Below are the generalized Jacobi
identities  for $n=0, 1, 2, 3$:
\begin{gather}
    \{\F\}=0,\\
    \{\{a\}\}+\{\F,a\}=0,  \label{eqjac1}\\
    \{\{a,b\}\}+\{\{a\},b\}+(-1)^{\at\bt}\{\{b\},a\}+\{\F,a,b\}=0,
    \label{eqjac2}
    \displaybreak[0]\\
\begin{split} \label{eqjac3}
    \{\{a,b,c\}\}
   +  \{\{a,b\},c\}+(-1)^{\bt\ct}\{\{a,c\},b\}+(-1)^{\at(\bt+\ct)}\{\{b,c\},a\}
    \\
   \, + \{\{a\},b,c\}+(-1)^{\at\bt}\{\{b\},a,c\}+(-1)^{(\at+\bt)\ct}\{\{c\},a,b\} \\
   +  \{\F,a,b,c\}=0.
\end{split}
\end{gather}

We shall call the $L_{\infty}$-algebras with $\F=0$, \textit{strict}.

For strict $L_{\infty}$-algebras, the Jacobi identities start from
$n=1$, and in~\eqref{eqjacn} the summation is over $k\geq 0,\, l> 0$.
The identities~\eqref{eqjac1}--\eqref{eqjac3} for strict
$L_{\infty}$-algebras simplify to
\begin{gather}
    d^2a  =0,  \label{eqjac1strict}\\
    d\{a,b\} +\{da ,b\}+(-1)^{\at\bt}\{db ,a\}=0,
    \label{eqjac2strict}\\
\begin{split} \label{eqjac3strict}
    d\{a,b,c\}
   +  \{\{a,b\},c\}+(-1)^{\bt\ct}\{\{a,c\},b\}+(-1)^{\at(\bt+\ct)}\{\{b,c\},a\}
    \\
    +\, \{da ,b,c\}+(-1)^{\at\bt}\{db ,a,c\}+(-1)^{(\at+\bt)\ct}\{dc,a,b\} =0,
\end{split}
\end{gather}
if we denote the unary bracket as $d:=\{\_\}$. That is, $d$ acts as a
differential, it has  the derivation property w.r.t. the binary
bracket, and the usual Jacobi holds for the binary bracket with a
homotopy correction. The identities with $n>3$ impose extra relations
for this homotopy and all the higher homotopies (hence `strongly'  in
the name).

As the operations $\{a_1,\ldots,a_n\}$ are multilinear and symmetric,
they are completely determined by the values on coinciding even
arguments: $\{\x,\ldots,\x\}$ where $\x$ is an even element of $V$
(to this end,  extension of scalars by odd constants should be
allowed). A generating function for these operations can be
conveniently written as a (formal) odd vector field on the vector
space $V$ considered as a supermanifold:
\begin{equation}\label{eqgenfunforbracket}
    Q=Q^i(\x)\der{}{\x^i}:=
    \sum_{n\geq 0}\frac{1}{n!}\,\{\underbrace{\x,\ldots,\x}_{n}\}.
\end{equation}
The elements of $V$ are identified with (constant) vector fields  as
$u=u^ie_i\,\leftrightarrow \, u^i \p_i$. If we denote the
\textit{$n$-th Jacobiator}, i.e., the \LHS  of~\eqref{eqjacn}, by
$J^n(a_1,\ldots,a_n)$,  it is clear that $J^n$ also give multilinear
symmetric operations on $V$. Hence they are, too, defined by their
values on equal even arguments.  The expression simplifies greatly,
and we have
\begin{equation}\label{eqjacobiator}
    J^n(\x,\ldots,\x)=\sum_{l=0}^{n}\frac{n!}{l!\,(n-l)!}\,
     \bigl\{\{\underbrace{\x,\ldots,\x}_{n-l}\},\underbrace{\x,\ldots,\x}_{l}\bigr\}
\end{equation}
for an even $\x$.  Abbreviating $J^n(\x,\ldots,\x)$ to $J^n(\x)$
we can write a generating function as
\begin{equation}
    J:=\sum_{n\geq 0} \frac{1}{n!}\,J^n(\x),
\end{equation}
which is an even (formal) vector field on the supermanifold $V$. One
can directly see that $J=Q^2=\frac{1}{2}\,[Q,Q]$. Hence all the
Jacobi identities can be compactly written as $Q^2=0$. Notice that
for strict $L_{\infty}$-algebras the vector fields  $Q$ and $J=Q^2$
vanish at the origin.

\begin{rem} \label{remdefshla}
There is a difference between the sign conventions of our
Definition~\ref{defshL} and the `standard' definitions of the
$L_{\infty}$-algebras as in~\cite{lada:stasheff, lada:shla}. It comes
from two sources. First, there is a choice between the `graded' (=
$\ZZ$-graded) and `super' viewpoints. Second, in supermathematics one
can choose between `symmetric' and `antisymmetric' constructions
using the parity shift. In~\cite{lada:stasheff, lada:shla} all vector
spaces are $\ZZ$-graded, but not `super', and brackets are
antisymmetric in the graded sense, i.e., involving the usual signs of
permutations together with the `Koszul signs' coming from the
$\ZZ$-grading. We prefer to work with the `super' conventions where
all the signs come from the $\Z$-grading (but not from any extra
$\ZZ$-grading be it present), and our brackets are (super) symmetric.
This has an advantage that it allows to use geometric language and
certain signs are simplified (e.g., the signs of permutations do not
enter). On the other hand, the definitions in~\cite{lada:stasheff,
lada:shla}  include directly the ordinary Lie algebras as a
particular case. A passage from~\cite{lada:stasheff, lada:shla} to
our conventions consists in introducing a $\Z$-grading (parity) as
the degree mod $2$ and applying the parity shift. Notice that it
reverses the parities of brackets with even numbers of arguments and
turns antisymmetric operations into symmetric. More precisely, let
$\Pi$ be the parity reversion functor. Suppose  $V=\Pi \mathfrak{g}$.
If we relate operations in $V$ and $\mathfrak{g}$ by the equality
\begin{equation*}
    \P [x_1,\ldots,x_n]=  \{\P x_1,\ldots,\P x_n\}\,(-1)^{\e},
\end{equation*}
$x_i\in \mathfrak{g}$, where $\e=\xt_1(n-1)+\dots+\xt_{n-1}$, then
(assuming that all brackets in $V$ are odd), the brackets in
$\mathfrak{g}$ with an even number of arguments will be even and with
odd will be odd; the antisymmetry of brackets in $\mathfrak{g}$ is
equivalent to the symmetry of brackets in $V$; the Jacobi identities
in the form of~\cite{lada:stasheff, lada:shla}  for the brackets in
$\mathfrak{g}$, extending the ordinary Jacobi identity for Lie
algebras, are equivalent to the Jacobi identities in the
form~\eqref{eqjacn} for the brackets in $V$. One might prefer to call
such a $V=\P \mathfrak{g}$, an `$L_{\infty}$-antialgebra'. However,
we  shall stick to Definition~\ref{defshL} throughout this paper.
Notice that our conventions are close to those
in~\cite{zwiebach:93csft}.
\end{rem}

\begin{rem}
In almost all standard approaches to $L_{\infty}$-algebras there is
no $0$-ary bracket or, rather, it is assumed that the corresponding
element $\F=\{\void\}$ is zero. (Except in~\cite{zwiebach:93csft} and
some other physical works; the algebras with a non-zero
$\F=\{\void\}$  are called sometimes `weak' or `with background'.)
Hence the standard $L_{\infty}$-algebras are always `strict' in our
sense. We have allowed for a $0$-ary operation because $\F\neq 0$
does occur naturally in some  our examples, and even where it does
not, including it  sometimes simplifies the exposition.
\end{rem}

\section{Main theorem} \label{secmain}

Let us return to the setting described above, i.e., a Lie
superalgebra $L$ with a projector $P$ on an Abelian subalgebra, and
an element $\D\in L$. Forget for a moment about restrictions on $\D$.

\begin{de}
For an arbitrary element $\D\in L$, even or odd, we call the
\textit{$n$-th derived bracket} of $\D$  the following operation on
the subspace $V:=P(L)\subset L$:
\begin{equation}\label{eqderived}
    \{a_1,\ldots,a_n\}_{\D}:=P[\ldots[[\D,a_1],a_2],\ldots,a_n]\,,
\end{equation}
where $a_i\in V$. Here $n=0,1,2,3,\ldots, \ $.
\end{de}

We get a set of $n$-ary operations~\eqref{eqderived} on the space
$V$.   Clearly, they are multilinear and of the same parity as
$\D$. (For $n=0$,  we get $\F:=\{\void\}=P(\D)$.) Notice that they
are always symmetric. Indeed, for the interchange of $a_1$ and
$a_2$, since $[a_1,a_2]=0$, we have
$[[\D,a_1],a_2]=-(-1)^{\e\at_1}[a_1,[\D,a_2]]=(-1)^{\e\at_1+\at_1(\e+\at_2)}[[\D,a_2],a_1]=
(-1)^{\at_1\at_2}[[\D,a_2],a_1]$, where $\e=\tilde \D$. Hence
\begin{equation*}
    \{a_1,a_2,\ldots,a_n\}_{\D}=(-1)^{\at_1\at_2}\{a_2,a_1,\ldots,a_n\}_{\D}.
\end{equation*}
Similarly for other adjacent arguments. For the coinciding even
arguments of the $n$-th derived bracket we have
\begin{equation}\label{eqderived2}
    \{\underbrace{\x,\ldots,\x}_{n}\}_{\D}=P (-\ad\x)^n\D
\end{equation}
(which is reminiscent of the $n$-th derivative of an 
$f(x)$ at a point $x_0$).

In the sequel we shall be particularly interested in the derived
brackets of an odd element $\D$ and  of its square
$\D^2=\frac{1}{2}[\D,\D]$.

Let $\D\in L$ be odd. Consider the Jacobiators $J^n_{\D}(\x)$ for the
derived brackets of $\D$. From~\eqref{eqjacobiator}
and~\eqref{eqderived2} we get
\begin{multline}\label{eqjacobforderived}
    J^n_{\D}(\x)=
    (-1)^n\sum_{l=0}^{n}\frac{n!}{l!\,(n-l)!}\,P(\ad\x)^l[\D,P(\ad\x)^{n-l}\D]=\\
    (-1)^n\sum_{l=0}^{n}\frac{n!}{l!\,(n-l)!}\,P[(\ad\x)^l\D,P(\ad\x)^{n-l}\D],
\end{multline}
where to obtain the second equality we used the Leibniz formula for
$(\ad \x)^l$ w.r.t. the Lie bracket in $L$ and the vanishing of the
commutators between elements of $V\subset L$.

\begin{thm}\label{thmmain}
Suppose $P$  satisfies~\eqref{eqabelian}, \eqref{eqdistrib}. Let
$\D$ be an arbitrary odd element. Then the $n$-th Jacobiator
$J^n_{\D}$ for the derived brackets of  $\D$ is exactly the $n$-th
derived bracket of $\D^2$:
\begin{equation}
    J^n_{\D}(a_1,\ldots,a_n)=\{a_1,\ldots,a_n\}_{\D^2}.
\end{equation}
\end{thm}
\begin{proof}
We shall prove the required identity for the coinciding even
arguments:
\begin{equation}
    J^n_{\D}(\x)=\{\underbrace{\x,\ldots,\x}_{n}\}_{\D^2}.
\end{equation}
Indeed, for the \LHS we can apply
formula~\eqref{eqjacobforderived}.  Let us analyze the cases of
$n$ odd and $n$ even separately. Suppose $n=2m+1$. Then we have
\begin{multline*}
    -J^{2m+1}_{\D}(\x)=P[\D,P(\ad\x)^{2m+1}\D]+
    \frac{(2m+1)!}{1!\,(2m)!}\,P[\ad\x.\D,P(\ad\x)^{2m}\D]+ \\
     \ldots+
    \frac{(2m+1)!}{m!\,(m+1)!}\,P[(\ad\x)^m\D,P(\ad\x)^{m+1}\D]+
    \\
    \frac{(2m+1)!}{(m+1)!\,m!}\,P[(\ad\x)^{m+1}\D,P(\ad\x)^{m}\D]+ \\
    \ldots+
    \frac{(2m+1)!}{(2m)!\,1!}\,P[(\ad\x)^{2m}\D,P(\ad\x.\D)]+ P[(\ad\x)^{2m+1}\D,P(\D)].
\end{multline*}
The terms corresponding to $l$ and $2m+1-l$, where
$l=0,1,\ldots,m$ can be grouped in pairs, and to each of the pairs
we can apply the distributive law~\eqref{eqdistrib}. Thus we get
after taking  $P$ out:
\begin{multline*}
    J^{2m+1}_{\D}(\x)=-P\Biggl(\sum_{l=0}^{m}\frac{(2m+1)!}{l!\,(2m+1-l)!}\, [(\ad\x)^l\D, (\ad\x)^{2m+1-l}\D]\Biggr)
    =\\
    -\frac{1}{2}\,P\Biggl(\sum_{l=0}^{2m+1}\frac{(2m+1)!}{l!\,(2m+1-l)!}\, [(\ad\x)^l\D, (\ad\x)^{2m+1-l}\D]\Biggr)
    =\\
    -\frac{1}{2}\,P(\ad\x)^{2m+1}[\D,\D]=P(-\ad\x)^{2m+1}\D^2.
\end{multline*}
Here we used the Leibniz identity for $(\ad\x)^{2m+1}$ w.r.t. the
commutator in $L$. Now suppose $n=2m>0$. We have
\begin{multline*}
    + J^{2m}_{\D}(\x)=P[\D,P(\ad\x)^{2m}\D]+
    \frac{(2m)!}{1!\,(2m-1)!}\,P[\ad\x.\D,P(\ad\x)^{2m-1}\D]+ \\
     \ldots+
    \frac{(2m)!}{(m-1)!\,(m+1)!}\,P[(\ad\x)^{m-1}\D,P(\ad\x)^{m+1}\D]+
    \\
    \frac{(2m)!}{m!\, m!}\,P[(\ad\x)^m\D,P(\ad\x)^{m}\D]+
    \\
    \frac{(2m)!}{(m+1)!\,(m-1)!}\,P[(\ad\x)^{m+1}\D,P(\ad\x)^{m-1}\D]+  \ldots\\
   +
    \frac{(2m)!}{(2m-1)!\,1!}\,P[(\ad\x)^{2m-1}\D,P(\ad\x.\D)]+
    \frac{(2m)!}{(2m)!\,0!}\,P[(\ad\x)^{2m}\D,P(\D)].
\end{multline*}
All  terms except for the term with $l=m$ can be grouped in pairs
and transformed as above. To the term corresponding to $l=m$ we
can apply the identity $P[a,a]=2P[Pa,a]$, which follows from the
distributive law~\eqref{eqdistrib}, valid for any odd $a\in L$.
Hence we get, similarly to the above,
\begin{multline*}
    J^{2m}_{\D}(\x)=
    P\Biggl(\sum_{l=0}^{m-1}\frac{(2m)!}{l!\,(2m-l)!}\, [(\ad\x)^l\D,
    (\ad\x)^{2m-l}\D]\Biggr)+ \\
    \frac{1}{2}\,\frac{(2m)!}{m!\, m!}\,P[(\ad\x)^m\D,
    (\ad\x)^{m}\D]= \\
    \frac{1}{2}\, P\Biggl(\sum_{l=0}^{2m}\frac{(2m)!}{l!\,(2m-l)!}\, [(\ad\x)^l\D,
    (\ad\x)^{2m-l}\D]\Biggr)= \\
    \frac{1}{2}\,P(\ad\x)^{2m}[\D,\D]=P(\ad\x)^{2m}\D^2.
\end{multline*}
For completeness, notice that for $n=0$ we have
$J^0_{\D}=\{\{\void\}_{\D}\}_{\D}=P[\D,P(\D)]=\frac{1}{2}\,P[\D,\D]=P(\D^2)=\{\void\}_{\D^2}$.
We conclude that in all cases
\begin{equation}
    J^n_{\D}(\x)=P(-\ad\x)^{n}\D^2=\{\underbrace{\x,\ldots,\x}_{n}\}_{\D^2},
\end{equation}
as claimed.
\end{proof}

\begin{cor} \label{corlinfty}
In the setup of Theorem~\ref{thmmain}, if $\D^2=0$, the derived
brackets of $\D$ make $V$  an $L_{\infty}$-algebra.
\end{cor}

This allows a  generalization, which naturally comes up in
examples.

\begin{de}
For any element $\D\in L$ we define the number $r$ to be the
\textit{order of $\D$ w.r.t. a subalgebra $V\subset L$} if all
$(r+1)$-fold commutators $[\ldots[\D,a_1],\ldots,a_{r+1}]$ with
arbitrary elements of $V$ identically vanish. Notation: $\ord_V
\D$.
\end{de}

  This is a filtration in $L$.

\begin{cor}\label{corhigherlinfty}
In the setup of Theorem~\ref{thmmain}, if  $\,\ord_V \D^2\leq
r\,$,   then the derived brackets of $\D$ satisfy the Jacobi
identities of orders $n>r$.
\end{cor}

We call the algebras given by Corollary~\ref{corhigherlinfty},
\textit{$L_{\infty}$-algebras of order $>r$}.

Notice that any higher Jacobi identity includes all $n$-brackets with
$n=0,1,\ldots$ . As above, we can speak about \textit{strict}
$L_{\infty}$-algebras of order $>r$ if the $0$-bracket $\F$ vanishes.
A natural question is, when one can split the element $\F=P(\D)$ from
the Jacobi identities of orders $n\geq 1$ and simply drop the $0$-ary
bracket from consideration. This happens if $\F$ is an annihilator of
all $n$-brackets, $n=2,3,\ldots$ . Besides an evident case
$\F=P(\D)=0$, a sufficient condition is $P[\D,P(\D)]=[\D,P(\D)]$. See
examples in the next section.

\section{Applications}

\label{secappl}

In this section, we consider some examples of applications of
Theorem~\ref{thmmain} and Corollaries~\ref{corlinfty},
\ref{corhigherlinfty}.

\begin{ex}\label{exhomfield}
Consider a vector space $V$, with the algebra  $L=\Vect (V)$ and  the
projector $P$ as in Example~\ref{exvectorfields}.  Take as $\D$ an
arbitrary odd vector field $Q\in \Vect (V)$,
\begin{equation*}
Q=Q^k(\x)\,\der{}{\x^k}=\left(Q_0^k+\x^i
\,Q_i^k+\frac{1}{2}\,\x^j\x^i\, Q_{ij}^k+
\frac{1}{3!}\,\x^l\x^j\x^i\, Q_{ijl}^k+\ldots\right)\der{}{\x^k}\,.
\end{equation*}
The derived brackets of $Q$,
\begin{equation} \label{eqqbracket}
\{a_1,\ldots,  a_n\}_Q= [\ldots[[Q,a_1],a_2],\ldots,a_n](0)\,,
\end{equation}
where $a_i\in V$ are identified with the corresponding constant
vector fields, are given by the coefficients of the Maclaurin
expansion:
\begin{gather*}
Q_0=Q_0^k e_k, \\ 
de_i:=\{e_i\}= (-1)^{\itt+1}\,Q_i^k e_k, \quad
\{e_i,e_j\}=(-1)^{\itt+\jtt}\,Q_{ij}^k e_k, \\
\{e_i,e_j,e_l\}=(-1)^{\itt+\jtt+\lt+1}\,Q_{ijl}^k e_k, \  \ldots \ ,
\end{gather*}
for the basis $e_i$. Here we denoted $\itt=\tilde e_i$. These  are
precisely the brackets on $V$ for which the vector field $Q$ (or,
rather, its Maclaurin series) is the generating
function~\eqref{eqgenfunforbracket}. Hence for $Q^2=0$ we recover the
$1-1$-correspondence between $L_{\infty}$-algebras and homological
vector fields. Moreover, we see that it is given by the explicit
formula~\eqref{eqqbracket}. Set $Q_0=0$. Then the algebra is strict.
For vector fields on $V$, the order w.r.t. the subalgebra of constant
vector fields $V\subset \Vect(V)$ is the degree in the variables
$\x^i$ (as a filtration). We conclude that strict
$L_{\infty}$-algebras of order $>r$ are in a $1-1$-correspondence
with odd vector fields $Q$ vanishing at the origin with the square
$Q^2$ of degree $\leq r$ in coordinates.
\end{ex}

\begin{ex} \label{exham}
For a (super)manifold $M$, consider   $\fun(M)\subset \fun(T^*M)$ as
in Example~\ref{exhamiltonians}. The projector is the pullback. Any
odd Hamiltonian $S\in \fun (T^*M)$ defines a sequence of
\textit{higher Schouten \textnormal{(}= odd Poisson\textnormal{)}
brackets} in $\fun(M)$ by the formula
\begin{equation*}
    \{f_1,\ldots,f_n\}_S:=  (\ldots((S,f_1),f_2),\ldots,f_n)_{|p=0}
\end{equation*}
(the parentheses stand for the canonical Poisson bracket on $T^*M$).
Here $f_1,\ldots,f_n\in\fun(M)$. They satisfy the Jacobi identities
of all orders if $(S,S)=0$. Notice that a Hamiltonian has a finite
order w.r.t. the subalgebra $\fun(M)$ if it is polynomial in  $p_a$,
and the order is the respective degree. The Jacobi identities can be
obtained one by one by putting restrictions on the order of $(S,S)$.
If
\begin{equation*}
S=S(x,p)=S_0(x)+ S^a(x) p_a +\frac{1}{2}\, S^{ab}(x) p_bp_a+
\frac{1}{3!}\,S^{abc}(x) p_cp_bp_a+\ldots\,,
\end{equation*}
then
\begin{gather*}
\{\void\}=S_0, \quad \delta f:=\{f\}= S^a\,\p_a f, \quad
\{f,g\}=(-1)^{\ft\at}\,S^{ab} \,\p_bf\,\p_ag, \\
\{f,g,h\}= (-1)^{\ft(\at+\bt)+\gt\at}\,S^{abc}
\,\p_cf\,\p_bg\,\p_a h, \ \ldots \ .
\end{gather*}
If $(S,S)$ is of degree $\leq r$ in $p_a$, then the brackets satisfy
the Jacobi identities of orders $\geq r+1$. In this example each of
the higher Schouten brackets is a multi-derivation, i.e., satisfies
the Leibniz rule w.r.t. the usual product, in each argument. Hence
the algebras that we obtain  are particular homotopy analogs of odd
Poisson (= Schouten, Gerstenhaber) algebras.  The `strict' case is
when $S_{|p=0}=0$.
\end{ex}

\begin{ex}\label{exmultivec}
Similarly to the above, take as $L$ the algebra of multivector
fields $\fun(\P T^*M)$ with the canonical Schouten bracket. Here
we have to change parity to obtain a Lie superalgebra. The rest
goes as in Example~\ref{exham}. Any even multivector field $P\in
\fun(\P T^*M)$ provides a sequence of  \textit{higher Poisson
brackets} in $\fun(M)$:
\begin{equation*}
    \{f_1,\ldots,f_n\}_P :=
    \lsch\ldots\lsch\lsch
    P,f_1\rsch,f_2\rsch,\ldots,f_n\rsch_{|x^*=0}
\end{equation*}
The brackets of odd orders are odd, the brackets of  even orders are
even. We have
\begin{equation*}
    \{f_1,\ldots,f_n\}_P =
    P^{a_1\ldots a_n}(x)\,\p_{a_n}f_n\ldots
    \p_{a_1}f_1
\end{equation*}
for even functions (for arbitrary functions the formula  follows by
linearity, using multiplication by odd constants), where
\begin{equation*}
P=P(x,x^*)=P_0(x)+P^a(x)\, x^*_a +\frac{1}{2}\, P^{ab}(x)\,
x^*_bx^*_a+ \frac{1}{3!}\,P^{abc}(x)\, x^*_cx^*_bx^*_a+\ldots\,,
\end{equation*}
with the full set of the Jacobi identities being equivalent to $\lsch
P,P\rsch=0$. Again, there is a possibility of getting the Jacobi
identities step by step by putting restrictions  on the degree of
$\lsch P,P\rsch$. As in Example~\ref{exham}, each of the higher
brackets strictly satisfies the Leibniz rule w.r.t. the product of
functions.
\end{ex}

Examples~\ref{exmultivec} and~\ref{exham} generalize  classical
Poisson and Schouten (= odd Poisson) structures,  as
Example~\ref{exhomfield}  generalizes classical Lie algebras. Indeed,
for a bivector field $P=\frac{1}{2}\,P^{ab}x^*_bx^*_a$ or  an odd
Hamiltonian quadratic in the momenta $S=\frac{1}{2}\, S^{ab} p_bp_a$,
the binary derived  bracket is an ordinary Poisson or Schouten
bracket, respectively, and all other brackets vanish. (Similarly,
after the shift of parity the bracket in a Lie algebra is the binary
derived bracket for a homological vector field
$Q=\frac{1}{2}\,\x^i\x^j\,Q^k_{ji}\,\p_k$ that is quadratic in
coordinates.) A mechanism for the arising of higher brackets can be
to take a quadratic Hamiltonian or a bivector field    generating an
ordinary Schouten or Poisson  bracket, and apply to it  a canonical
transformation that fixes the zero section but   not   the bundle
structure.

Notice in both examples the possibility of obtaining higher odd or
even Poisson brackets from a non-polynomial Hamiltonian or
multivector field (the latter is possible only in the super case). It
is the Taylor expansion around the zero section $M\subset T^*M$ or
$M\subset \Pi T^*M$ that counts.

\begin{ex}\label{exdelta}
Consider a  commutative associative algebra with a unit $A$, e.g.,
an algebra of smooth functions $\fun(M)$. Let $L$ be the algebra
of all linear operators in $A$ w.r.t. the   commutator, and $V=A$
considered as an Abelian subalgebra in $L$. Let $P\co \End A\to
\End A$ be the evaluation at $1$, as in Example~\ref{exoperators}.
Let $\D$ be an arbitrary odd operator in $A$.  The derived
brackets of $\D$,
\begin{equation} \label{eqdbracket}
 \{f_1,\ldots,f_n\}_{\D}:=
 [\ldots[[\D,f_1],f_2],\ldots,f_n](1),
\end{equation}
will be, respectively,
\begin{align*}
    \{\void \}&=\F =\D 1, \\
    \{f\}_{\D}&=\D'f=\D f-\D 1 \cdot f,\\
    \{f,g\}_{\D}&=\D(fg)- \D f \cdot g-(-1)^{\ft} f\cdot  \D g +\D 1 \cdot fg, \\
    \{f,g,h\}_{\D}&=
    \D(fgh)-\D(fg)\cdot h
    -(-1)^{\gt\hht}\D(fh)\cdot g \\
        &  \quad   -(-1)^{\ft(\gt+\hht)}\D(gh)\cdot f +
    \D f\cdot gh  \\
    &  \quad + (-1)^{\ft\gt}\D g\cdot fh+(-1)^{\hht(\ft+\gt)}\D h\cdot fg-
    \D1\cdot fgh \\
    \ldots & \ldots   \ldots
\end{align*}
One can check that these brackets satisfy the following identity
w.r.t. the product of functions:
\begin{multline*}
    \{f_1,\ldots,f_{n-1},gh\}_{\D}=\\
    \{f_1,\ldots,f_{n-1},g\}_{\D}h+(-1)^{\gt\hht}
    \{f_1,\ldots,f_{n-1},h\}_{\D}g +
    \{f_1,\ldots,f_{n-1},g,h\}_{\D}\,,
\end{multline*}
i.e., the $(n+1)$-th bracket arises as the failure of the Leibniz
rule for the $n$-th bracket.  If $\D$ is a differential operator of
order $s$, then the $(s+1)$-th bracket and all higher brackets
identically vanish, and the $s$-th bracket is a (symmetric)
multi-derivation of the algebra $A$. (It is nothing but the
polarization of the principal symbol of $\D$.) The usual order of a
differential operator is exactly the `order w.r.t. the subalgebra
$A$'. The $k$-th bracket with $1\leq k\leq s$ is in this case a
differential operator of order $s-k+1$ in each argument. One can view
these brackets  as  consecutive ``polarizations'' of the operator
$\D$. It is instructive to write them down explicitly for a
particular operator $\D$ in a differential-geometric setting (see
below). As follows from Theorem~\ref{thmmain}, if $\D^2=0$, then the
derived brackets of $\D$ satisfy the Jacobi identities of all orders;
otherwise, by requiring $\ord \D^2\leq r$ we obtain the Jacobi
identities of orders $r+1$ and higher.
\end{ex}

\begin{rem}
That the brackets~\eqref{eqdbracket} give an $L_{\infty}$-algebra
if $\D^2=0$ was for the first time proved in~\cite{bering:higher},
by rather hard calculations.
\end{rem}

The $n$-brackets~\eqref{eqdbracket} with $n\geq 1$ will not change if
we replace $\D$ by $\D'=\D-\D1$.   Let $J'^n_{\D}$ denote the $n$-th
Jacobiator with $\F$ dropped, and $J^n_{\D}$ stands for the full
Jacobiator, $n>0$. Then $J'^n_{\D}=J^n_{\D'}$. Applying
Theorem~\ref{thmmain}, we identify $J^n_{\D}$ and $J'^n_{\D}$ with
the $n$-brackets generated by $\D^2$ and $\D'^2$ respectively. Since
$\D'^2=\D^2-[\D,\D 1]$, by comparing the orders we conclude that
$J'^n_{\D}=J^n_{\D}$ for all $n=s, s+1,\ldots,2s-1$ if $\ord \D\leq
s$. Hence $\F=\D 1$ can be dropped from the $n$-th Jacobi identity
for the brackets generated by $\D$ exactly for these numbers $n$.

The construction in Example~\ref{exdelta}  generalizes  the
interpretation of a classical odd Poisson bracket as the derived
bracket of a `generating operator' of second order  (= an odd
Laplacian, a `BV-operator'). This approach   was particularly useful
for the analysis of second order differential operators
in~\cite{tv:laplace1,tv:laplace2} (see also \cite{tv:laplace2bis}).

\begin{ex}[\cite{tv:laplace1,tv:laplace2}] \label{exdeltasec}
If $\D$ is an odd $2$-nd order differential operator in $\fun(M)$, in
local coordinates
\begin{equation*}
    \D= R(x) + T^a(x)\,\p_a + \frac{1}{2}\, S^{ab}(x) \,\p_b\p_a\,,
\end{equation*}
then we get
\begin{align*}
    \F&=\D 1=R \\
    \{f\}_{\D}&=\D'=T^a \,\p_af + \frac{1}{2}\, S^{ab}  \,\p_b\p_af\,,\\
    \{f,g\}_{\D}&=(-1)^{\ft\at}\,S^{ab}\,\p_bf\,\p_ag\,.
\end{align*}
All the higher brackets vanish. Automatically $\ord \D^2\leq 3$. If
$\ord \D^2\leq 2$, then $\{f,g\}_{\D}$ satisfies the usual Jacobi
identity, making $\fun(M)$ into an odd Poisson algebra. If $\ord
\D^2\leq 1$, then $\D'$ is a derivation of the bracket. Finally, if
$\ord \D^2\leq 0$ {and} $\D 1=0$, then $\D=\D'$ is a differential;
the resulting algebraic structure  is known as a
\textit{Batalin--Vilkovisky algebra}. (Notice that $\D1$ does not
affect the Jacobi identities with $n=2, 3$.)
\end{ex}

\begin{ex}\label{exdeltathird}
For an odd $3$-rd order differential operator in $\fun(M)$, in local
coordinates
\begin{equation*}
    \D= R(x) + T^a(x)\,\p_a + \frac{1}{2}\, U^{ab}(x) \,\p_b\p_a+
    \frac{1}{3!}\,S^{abc}(x)\,\p_c\p_b\p_a\,,
\end{equation*}
we get
\begin{align*}
    \F&=\D 1=R \\
    \{f\}_{\D}&= \D'=T^a \,\p_af + \frac{1}{2}\, U^{ab}  \,\p_b\p_af+
    \frac{1}{3!}\,S^{abc} \,\p_c\p_b\p_af,\\
    \{f,g\}_{\D}&= (-1)^{\ft\at}\left( U^{ab}\,\p_bf\, \p_ag +
    \frac{1}{2}\,S^{abc}\bigl((-1)^{\ft\bt}\,\p_cf\, \p_b\p_ag+
    \p_c\p_bf\, \p_ag\bigr)\right)\,, \\
    \{f,g,h\}_{\D}&=(-1)^{\ft(\at+\bt)+\gt\at}\,S^{abc}\,\p_cf\,\p_bg\,\p_ah\,,
\end{align*}
and all the higher brackets vanish. Automatically $\ord \D^2\leq 5$.
Not affected by $\D 1$ are the Jacobi identities with $n=3, 4, 5$. If
$\ord \D^2\leq 4$, then there holds the $5$-th order Jacobi identity
\begin{equation*}
    \sum_\text{shuffles} \pm\,\{\{f,g,h\}_{\D},e,k\}_{\D}=0.
\end{equation*}
It involves only the ternary bracket. If $\ord \D^2\leq 3$, then
also holds  the $4$-th order Jacobi identity
\begin{equation*}
    \sum_\text{shuffles} \pm\,\{\{f,g,h\}_{\D},e\}_{\D}+
    \sum_\text{shuffles} \pm\,\{\{f,g\}_{\D},h,e\}_{\D}=0.
\end{equation*}
If $\ord \D^2\leq 2$,   then in addition holds the $3$-rd order
Jacobi identity:
\begin{multline*}
    \sum_\text{cycle} \pm\,\{\{f,g\}_{\D},h\}_{\D} \\
    \pm\,
    \D'\{f,g,h\}_{\D}
     \pm\,\{\D'f,g,h\}_{\D}
    \pm\,\{f,\D'g,h\}_{\D} \pm\,\{f,g,\D'h\}_{\D}=0.
\end{multline*}
If $\ord \D^2\leq 1$,   we  get the $2$-nd order Jacobi identity
involving $\D 1=R$, which now cannot be ignored:
\begin{equation*}
    \D'\{f,g\}_{\D}
     \pm\,\{\D'f,g\}_{\D}
    \pm\,\{f,\D'g\}_{\D}+\{\D 1, f, g\}=0.
\end{equation*}
Finally, if $\ord \D^2\leq 0$, we arrive at the $1$-st order
Jacobi identity in the form $(\D')^2f+\{\D 1,f\}=0$. We have to
impose  $\D 1=0$ to get strictness back.
\end{ex}

\begin{rem}
The algebraic structure consisting of  all higher derived brackets of
an odd  differential operator of order $n$ and the usual
multiplication, should be considered  an example of a
\textit{homotopy Batalin--Vilkovisky algebra} (see~\cite{aa:hga,
kravchenko:bv, tamarkintsygan:bv} and a discussion in
Section~\ref{secdiscuss}).
\end{rem}

The behaviour of the brackets  in Example~\ref{exham} and
Example~\ref{exdelta}  w.r.t. the multiplication, at the first glance
seems very different. However, the identities satisfied by the
algebras obtained in Example~\ref{exham} can be seen as the
``classical limit'' of the identities for the algebras obtained in
Example~\ref{exdelta}. Indeed, if we redefine the brackets in
Example~\ref{exdelta} by inserting  the ``Planck's constant''
$\hbar$, as
\begin{equation*}
 \{f_1,\ldots,f_n\}_{\D}:=(-i\hbar)^{-n}
 [\ldots[[\D,f_1],f_2],\ldots,f_n](1)\,,
\end{equation*}
then they will satisfy the same Jacobi-type identities as before, but
the ``Leibniz identity'' will now read
\begin{multline*}
    \{f_1,\ldots,f_{n-1},gh\}_{\D}=
    \{f_1,\ldots,f_{n-1},g\}_{\D}h+(-1)^{\gt\hht}
    \{f_1,\ldots,f_{n-1},h\}_{\D}g   \\
   + (-i\hbar)\{f_1,\ldots,f_{n-1},g,h\}_{\D}\,,
\end{multline*}
which clearly becomes the strict derivation property when $\hbar\to
0$.

\section{Case of non-inner derivations} \label{secnonin}

Higher derived brackets generated by an element $\D$ naturally arise
in applications, as we saw it in the previous section. However, from
theoretical reasons and from the viewpoint of further generalizations
it seems natural to look also into a possibility to obtain a similar
construction from an arbitrary derivation of the superalgebra $L$
rather than  inner derivations given by $\D\in L$. It is indeed
possible and in particular allows to look at higher derived brackets
from yet another angle. Here we shall briefly outline the
construction and statements, leaving a more detailed exposition for
another occasion.

As above, let $L$ be a Lie superalgebra and $P$ a projector
satisfying  the identities $[Pa,Pb]=0$ and
\begin{equation*}
P[a,b]=P[Pa,b]+P[a,Pb]
\end{equation*}
for all $a,b\in L$. Recall that it means that both subspaces $V=\Im
P$ and $K=\Ker P$  are subalgebras and $V$ is Abelian. Consider an
arbitrary, even or odd, derivation $d$ of the Lie superalgebra $L$.
Let us assume that the kernel $K$ of $P$ is closed under $d$; this is
equivalent to the identity
\begin{equation}\label{eqpdp}
    PdP=Pd\,.
\end{equation}
(Notice that we do not assume the image of $P$, i.e., the subspace
$V$, to be closed under $d$.)

\begin{de}
The \textit{$n$-th   derived bracket} of $d$ is the following
operation on the subspace $V\subset L$:
\begin{equation}\label{eqderivednonin}
    \{a_1,\ldots,a_n\}_{d}:=P[\ldots[da_1,a_2],\ldots,a_n] \,,
\end{equation}
where $a_i\in V$. Here $n=1,2,3,\ldots\ $.
\end{de}

\begin{Rem}
If $V$ happens to be closed under $d$, then all the
$n$-brackets~\eqref{eqderivednonin} with $n>1$, will vanish. So it is
the non-commutativity of $d$ with $P$ that is the source of higher
derived brackets.
\end{Rem}

Brackets~\eqref{eqderivednonin} are even or odd depending on the
parity of $d$. Notice that there is no  $0$-ary bracket, differently
from the construction based on  $\D\in L$. Exactly as above follows
(from the derivation property of $d$, the Jacobi identity in $L$ and
the condition that the subalgebra $V\subset L$ is Abelian) that all
higher brackets~\eqref{eqderivednonin} are symmetric in the
$\Z$-graded sense.

\begin{thm} \label{thmmain2} Suppose $d$ is an odd derivation. Then the
$n$-th Jacobiator
of the derived brackets of $d$ is exactly the $n$-th derived bracket
of $d^2$:
\begin{equation}
    J^n_{d}(a_1,\ldots,a_n)=\{a_1,\ldots,a_n\}_{d^2}.
\end{equation}
Here $n=1,2 , 3,\ldots \ $. \textnormal{(In the formula for the
Jacobiator the $0$-th bracket should be set to zero.)}
\end{thm}

In particular, if $d^2=0$, the higher derived  brackets of $d$ make
the subspace $V$ a strict $L_{\infty}$-algebra. Clearly, it is also
possible to weaken the condition $d^2=0$ by considering instead of it
a filtration by an `order' of operators w.r.t. the subspace $V$, as
we did above for $\D$.

Theorem~\ref{thmmain2} is a generalization of Theorem~\ref{thmmain}
if  $P(\D)=0$. As we have seen it in the examples, this not always
the case, so better to consider these two statements as independent,
though closely related.

The construction of higher derived brackets from an arbitrary
derivation makes it possible to give for them a nice
homotopy-theoretic interpretation, as follows\footnote{A
homotopy-theoretic interpretation of our original construction with
$\D$ was conjectured by an anonymous referee of this paper, who
proposed to extend  the brackets generated by  $\D$ by  formulae
similar to ~\eqref{eqextbracket}--\eqref{eqextbracketlast} deduced
below.}. Let $d$  be an odd derivation of the superalgebra Lie $L$
such that $d^2=0$. So $L$ together with $d$ is a differential Lie
superalgebra. Consider the subalgebra $K=\Ker P$ in $L$. It is a
differential subalgebra. Consider the inclusion map $i\co K\to L$.
Forget for a moment about the algebra structure and consider it just
as an inclusion of complexes. (For our purposes, a \textit{complex}
is a $\Z$-graded vector space with an odd endomorphism of square
zero.) As topologists know, ``every map can be made a fibration'', by
applying a cocylinder construction. The fiber of this fibration is
known as a `homotopy fiber' of the original map. What is a homotopy
fiber for the inclusion $i\co K\to L$? The claim is, it is the space
$\Pi V$. Moreover, the higher derived brackets will make it a
homotopy fiber in the category of $L_{\infty}$-algebras. More
precisely, the following statements hold.

Let $i\co K\to L$ be an arbitrary inclusion of complexes such that
there is given a complementary subspace $V\subset L$ for $K$, so that
$L=K\oplus V$. ({$V$ is not necessarily a subcomplex.}) Let $P$ be
the projector onto $V$ parallel to $K$. The space $V$ becomes a
complex with the differential $Pd$. Introduce into $L\oplus \Pi V$ an
operator $D$ as follows:
\begin{equation}\label{eqdincocyl}
    D(x,\Pi a):=\bigl(dx, -\Pi P(x+da)\bigr), 
\end{equation}
for $x\in L$, $a\in V$. Then $D^2=0$ (check!). Consider the maps
$j\co K\to L\oplus \Pi V$ and $p\co L\oplus \Pi V\to L$, where  $j\co
x\mapsto (x,0)$, $p\co (x,\Pi a)\mapsto x$.
\begin{lem}
The diagram
\begin{equation}\label{eqcocyl}
\begin{diagram}[small]
    K &       & \rTo^i          &       & L  \\
      & \rdTo_j &               & \ruTo_p &    \\
      &       & L\oplus \Pi V &       &    \\
\end{diagram}
\end{equation}
is a cocylinder diagram in the category of complexes, i.e., the maps
$j$ and $p$ are chain maps,   $i=p\circ j$, the map  $j\co K\to
L\oplus \Pi V$ is a monomorphism \textnormal{(`cofibration')} and a
quasi-isomorphism \textnormal{(`weak homotopy equivalence')}, and the
map $p\co L\oplus \Pi V\to L$ is an epimorphism
\textnormal{(`fibration')}.
\end{lem}

(A quasi-inverse for $j$ is $q\co (x,\Pi a)\mapsto (1-P)(x+da)$.)

It follows that $\Pi V=\Ker p$ taken with the differential $-\Pi Pd$
is a homotopy fiber or a co-cone of the inclusion of complexes $i\co
K\to L=K\oplus V$.

Now, if we come back to our original setup where $i\co K\to L$ is an
inclusion of differential Lie superalgebras, we want to provide the
cocylinder $L\oplus \Pi V$ with a bracket extending the one in $L$ so
that $j$ and $p$ will respect the brackets and $D$ 
be a derivation.  It turns out that this condition fixes the bracket
in $L\oplus \Pi V$ uniquely. In addition to the original Lie bracket
in $L$, appear new brackets between elements of $L$ and $\Pi V$, and
inside $\Pi V$:
\begin{equation*}
[x,\Pi a]:=(-1)^{\tilde x}\Pi P[x,a], \quad [\Pi a, \Pi
b]:=(-1)^{\tilde a+1}\Pi P[da,b]
\end{equation*}
Up to the parity shift, the latter bracket is immediately
recognizable as the beginning of our sequence of  higher derived
brackets generated by $d$ in $V$. The new binary bracket in $L\oplus
\Pi V$ does not  satisfy the Jacobi identity exactly; this gives rise
to ternary brackets $L\oplus \Pi V$ of the form similar to the above,
and so on. One can figure out the appearance of these higher brackets
by an incomplete induction. Since in this paper we work with
symmetric brackets, the final result is more conveniently formulated
after a parity shift. Applying $\Pi$ to~\eqref{eqcocyl} we get
\begin{equation}
\begin{diagram}[small]\label{eqcocylpi}
    \Pi K &       & \rTo^i          &       & \Pi L  \\
      & \rdTo_j &               & \ruTo_p &    \\
      &       & \Pi L\oplus V &       &    \\
\end{diagram}
\end{equation}
which is a cocylinder diagram for $i=i^{\Pi}\co \Pi K\to \Pi L$ in
the category of complexes. Here $D=D^{\Pi}$ in $\Pi L\oplus V$ is
$(\Pi x, a)\mapsto \bigl(-\Pi dx, P(x+da)\bigr)$. The desire to
extend the bracket in $\Pi L$ corresponding to the Lie bracket in $L$
keeping $D$ a derivation,  naturally leads to the following
definitions. The $0$-ary bracket in $\Pi L\oplus V$ is set to zero
and as the unary bracket we take the operator $D$:
\begin{align}\label{eqextbracket}
    \{\Pi x\}&=-\Pi dx+Px,\\
    \{a\}&=Pda\,.
\end{align}
Then we define the binary brackets as
\begin{align}
\{\Pi x, \Pi y\}&=\Pi[x,y](-1)^{\tilde x}, \\
 \{\Pi x, a\}&=P[x,a], \\
 \{a, b\}&= P[da,b].
\end{align}
The higher order brackets  we define as
\begin{align}
\{\Pi x, a_1,\ldots,a_n\}&=P[\ldots[x,a_1],\ldots,a_n], \\
\{a_1,\ldots,a_n\}&=P[\ldots[da_1,a_2],\ldots,a_n],
\label{eqextbracketlast}
\end{align}
where $n> 1$.  All the other brackets except obtainable from these by
symmetry, are set to zero. We arrive at a collection of odd symmetric
multilinear operations on $\Pi L\oplus V$. The subspace $V$ is an
ideal w.r.t. these operations and their restriction to $V$ coincides
with the higher derived brackets~\eqref{eqderivednonin}.
\begin{thm}\label{thmcocyl}
The operations~\eqref{eqextbracket}--\eqref{eqextbracketlast} make
the space $\Pi L\oplus V$ a strict $L_{\infty}$-algebra.
\end{thm}

It follows from Theorem~\ref{thmcocyl} that the diagram
~\eqref{eqcocylpi} is a cocylinder diagram in the category of
$L_{\infty}$-algebras, as it is clear that the maps $j$ and $p$
in~\eqref{eqcocylpi} strictly respect the brackets, in particular
giving $L_{\infty}$-maps.  As a corollary we see that $V$ considered
with the higher derived brackets~\eqref{eqderivednonin} is a homotopy
fiber of the inclusion of the (odd) differential Lie superalgebras
$\Pi K\to \Pi L$. If we change the viewpoint at $L_{\infty}$-algebras
and adopt a definition which differs from ours by the parity shift
(see Remark~\ref{remdefshla}), it will be possible to say that $\Pi
V$ (with the corresponding `shifted' higher derived brackets) is a
homotopy fiber for the inclusion $K\to L$.

The proofs of Theorems~\ref{thmmain2} and~\ref{thmcocyl}, and other
details, will be given elsewhere.

\section{Discussion}

\label{secdiscuss}

A derived bracket  (with this name) of two arguments appeared for the
first time in the paper by
Y.~Kosmann-Schwarzbach~\cite{yvette:derived}, who also
referred  to an unpublished text by Koszul of 1990. 
She proved that any odd derivation of a Loday (= Leibniz) algebra
generates a new Loday bracket of the opposite parity  by the formula
$[a,b]_D=(-1)^{\at}[Da,b]$. (The present author independently
introduced
a derived bracket 
around 1993  and proved a similar statement, in a slightly less
generality than~\cite{yvette:derived},  namely, without the Loday
algebras and working only with Lie superalgebras.) Unlike the
brackets introduced in the present paper, the bracket $[a,b]_D$ does
not necessarily satisfy (anti)symmetry even if the original bracket
does. Antisymmetry is restored on suitable subspaces or quotient
spaces provided  the derived bracket can be restricted there. The
present construction of higher derived brackets making use of a
projector $P$ solves the problem by forcing the bracket to remain in
a given subspace. The necessity to consider all the higher brackets,
not just the binary bracket, is the price.

Retrospectively,  binary derived brackets, considered on subspaces,
can be recognized in many constructions of differential geometry,
e.g., in the Cartan identities $[d,i_u]=\lie{u}$,
$[\lie{u},i_v]=i_{[u,v]}$ combined to give $i_{[u,v]}=[i_u,[d,i_v]]$.
An important example  is the coordinate-free expression for a Poisson
bracket generated by a bivector field $B$ via the canonical Schouten
bracket: $\{f,g\}_B=\lsch f,\lsch B, g\rsch\rsch$ up to a sign
depending on conventions, and a similar expression for a Schouten
structure via the canonical Poisson bracket. (For the author, these
expressions  were a starting point in the discovery of  derived
brackets.) Derived brackets have been also  used for describing Lie
algebroids (see, e.g.,~\cite{tv:graded}) and Courant
algebroids~\cite{roytenberg:courant, roytenberg:thesis}.

Higher brackets do not appear in these  classical  examples  because
for them the generating element is always, loosely speaking,
`quadratic': viz., a quadratic homological vector field
$Q=\frac{1}{2}\,\x^i\x^j\,Q^k_{ji}\p_k$, a bivector field
$B=\frac{1}{2}\,B^{ab}(x)x^*_bx^*_a$, a quadratic Hamiltonian
$S=\frac{1}{2}\,S^{ab}(x)p_bp_a$, an odd   Laplacian $\D$, etc. For
the same reason there is no need to introduce a projector to remain
in a chosen subspace of the `zero-order' elements (such as vector
fields with constant coefficients, functions on $M$ as opposed to
those on $T^*M$, zero-order operators, etc.). On the other hand, a
natural attempt to replace, say, a Poisson bivector field by an
arbitrary multivector field satisfying $\lsch P,P\rsch=0$ and still
have a bracket on functions, requires introducing a projector and
immediately leads to higher derived brackets and  homotopy analogs of
the classical examples. (See Example~\ref{exham} and other examples
in the previous section.)

The characterization of differential operators with the help of
multiple commutators  can be traced to
Grothendieck~\cite{groth:ega4}. Related to it the higher derived
brackets of Example~\ref{exdelta} essentially coincide with the
operations $\Phi^r$ introduced by Koszul
in~\cite{koszul:crochet85}.   For a differential operator $\D$ on
a graded commutative algebra, Koszul defined $\F^r_{\D}$ for $r>0$
as
\begin{equation*}
    \F^r_{\D}(a_1\otimes\ldots\otimes a_r):=m\circ (\D\otimes
    \id)\,\l^r (a_1\otimes\ldots\otimes a_r)
\end{equation*}
where  $\l^r (a_1\otimes\ldots\otimes a_r)=(a_1\otimes 1-1\otimes
a_1)\ldots (a_r\otimes 1-1\otimes a_r)$ and $m$ stands for the
multiplication; he stated that for each $r$, $\F^{r+1}_{\D}$ equals
the failure of the Leibniz identity for $\F^r_{\D}$. He was basically
interested in the binary operation $\F^2_{\D}$ generated by an odd
operator of second order. It was stated in~\cite{koszul:crochet85}
that the failure of the homotopy Jacobi identity for $\F^2_{\D}$
(involving $\F^3_{\D}$) equaled $\F^3_{\D^2}$, and that the Leibniz
identity for $\D$  and $\F^2_{\D}$ were equivalent to
$\F^2_{\D^2}=0$. Generalizations of Koszul's operations $\F^r$ for
various types of algebras, not necessarily commutative or associative
were studied  by F.~Akman~\cite{akman:genBV} (see
also~\cite{akman:masterid}). As in~\cite{koszul:crochet85}, she was
mainly concerned with the binary bracket $\F^2$.

``Higher antibrackets'' generated by an odd differential operator
$\D$, together with higher Poisson brackets, have appeared in the
series of physical papers~\cite{alfaro:nonabelian, bering:higher,
batalin:higher96, batalin:quanti98, batalin:open98,
batalin:dualities, batalin:general} motivated by a development of the
Batalin--Vilkovisky quantization. As Stasheff noted, they were also
hiding in works on Batalin--Fradkin--Vilkovisky formalism, such
as~\cite{kjeseth:brst}.  In~\cite{bering:higher} it was proved
directly that the higher brackets defined by~\eqref{eqdbracket} form
an $L_{\infty}$-algebra if $\D^2=0$. For ``general antibrackets'' on
differential operators defined in~\cite{batalin:general} as the
symmetrizations of multiple commutators
\begin{equation*}
[\ldots[[\D,A_1],A_2],\ldots,A_n]\,
\end{equation*}
(no evaluation at $1$, unlike~\eqref{eqdbracket}), where the
operators $A_i$ are arbitrary and do not have to belong to an Abelian
subalgebra, were obtained certain Jacobi-type identities more
complicated than those for the $L_{\infty}$-algebras. Such algebraic
structures are yet to be analyzed.

As we mentioned in Section~\ref{secappl}, higher derived brackets of
Example~\ref{exdelta} make the natural framework for the problem of
describing the generating operators of an odd bracket. Geometric
constructions related with these `Batalin--Vilkovisky operators' were
considered in~\cite{hov:deltabest}, \cite{ass:bv},
\cite{yvette:divergence},~\cite{hov:max, hov:semi, hov:proclms}. A
complete picture was obtained in~\cite{tv:laplace1, tv:laplace2}.
In~\cite{tv:laplace2}, H.~Khudaverdian and the author established a
one-to-one correspondence between second-order differential operators
on the algebra of densities $\V(M)$ on a supermanifold $M$ and binary
brackets in $\V(M)$.  For operators acting on functions, this
specializes to a correspondence between operators and pairs
consisting of a bracket on functions and an ``upper connection'' on
volume forms~\cite{tv:laplace2}. For  odd operators this gives a
description of the Jacobi  conditions in terms of this connection.
Constructions of the present paper, hopefully, can be useful for
generalizing the results of~\cite{tv:laplace2} to higher order
operators.

There were suggested different approaches  to Batalin--Vilkovisky
algebras ``up to homotopy'', as well as to homotopy Schouten (=
Gerstenhaber) algebras. Operadic approaches to the latter are
discussed in~\cite{aa:hga}. A  direct definition of a homotopy
Batalin--Vilkovisky algebra was suggested by
Olga~Kravchenko~\cite{kravchenko:bv}  and further generalized by
Tamarkin and Tsygan~\cite{tamarkintsygan:bv}. In particular, besides
the $L_{\infty}$-structure this definition provides for the (strong)
homotopy associativity of the product and  homotopy Leibniz
identities. The examples in Section~\ref{secappl} satisfy much
stricter conditions.  On the other hand, an example of higher
brackets of differential operators as in~\cite{batalin:general}
involves conditions that are weaker than those of an
$L_{\infty}$-structure. Therefore the final algebraic framework for
these notions is yet to be found.

The  higher derived brackets that we introduced here are not the most
general. A natural extension of our constructions should be to allow
the image of a projector $P$ to be an arbitrary Lie subalgebra, not
necessarily Abelian. A condition generalizing~\eqref{eqdistrib}
should then read
\begin{equation} \label{eqdistrib2}
    P[a,b]=P[Pa,b]+P[a,Pb]-[Pa,Pb].
\end{equation}
(Together with $P[Pa,Pb]=[Pa,Pb]$ that means that both  $\Im P$ and
$\Ker P$ are  subalgebras, i.e., the Lie superalgebra in question is
the sum of two subalgebras.) In such case  the symmetry of higher
derived brackets should remain only up to homotopy, and we should end
up with a yet more general  notion of a (strongly) homotopy Lie
algebra. In examples, this should lead also to more general cases of
homotopy Batalin--Vilkovisky algebras. For instance, when $P=\id$,
this should cover the ``general antibrackets''
of~\cite{batalin:general}. Projectors satisfying~\eqref{eqdistrib2}
appeared in~\cite{hov:campbell} with a totally different motivation.
It was shown, remarkably, that they come from   operators on an
associative algebra with a unit satisfying
\begin{equation}\label{eqhovik}
    P1=1, \quad P(a Pb)=P(ab), \quad P((Pa)b)=Pa\, Pb.
\end{equation}
Khudaverdian's result~\cite{hov:campbell} is that upon
conditions~\eqref{eqhovik}, the formal series
\begin{equation*}
\log (Pe^a)=\log \left(1+Pa+\frac{1}{2}\,P(a^2)+\ldots\right)
\end{equation*}
for any element $a$ of the associative algebra can be expressed via
commutators only and the action of $P$,  thus obtaining a
generalization of the Baker--Campbell--Hausdorff formula. (These
results were inspired by an analysis of certain Feynman diagrams in
quantum field theory. Examples, however, range to cobordism theory
and Novikov's operator doubles~\cite{hov:campbell}.) There must be a
connection with the present construction of higher derived brackets,
but it is yet to be understood. A question that is related, is to
give an analogous ``derived'' construction in the associative
setting, leading to $A_{\infty}$-algebras and their relatives.

Another interesting direction of  study should be derived brackets
and homotopy algebras arising from graded manifolds~\cite{tv:graded}
(see also~\cite{roytenberg:graded}).

We hope to  consider these questions elsewhere.



\def\cprime{$'$}

\end{document}